\documentclass[12pt]{amsart}
\usepackage{amsmath,amsfonts,amssymb,amsthm,amscd}
\usepackage{graphicx}
\usepackage{mathrsfs}
\usepackage{upgreek}
\usepackage[matrix,arrow,curve]{xy}
\usepackage{xcolor}

\usepackage{hyperref}
\definecolor{urlcolor}{HTML}{799B03} % цвет гиперссылок

\usepackage{bm}

\usepackage[T2A]{fontenc}
\usepackage[utf8]{inputenc}
\usepackage[english,russian]{babel}
\usepackage{bm}
\usepackage{mathrsfs}

\newcommand\hide[1]{\commented{gray}{Hidden:}{#1}}
\renewcommand\hide[1]\empty

\hide{
\unitlength=1mm 
\textwidth=390pt 
\textheight=550pt 
\topskip=10pt
\baselineskip=15pt
\topmargin=27pt 
\footskip=20pt
\oddsidemargin=39.5pt 
\marginparwidth=68pt 
\sloppy
\let\emptyset\varnothing
}

\newtheorem{thm}{Theorem}

\newtheorem{cor}{Corollary}

\theoremstyle{definition}
\newtheorem{df}{Definition}

\newtheorem{rem}{Remark}

\theoremstyle{remark}

\begin{document}

\renewcommand{\bm}{\boldsymbol}

\newcommand{\dom}{ {\mathop{\mathrm {dom}}\nolimits} }
\newcommand{\ran}{ {\mathop{\mathrm{ran}}\nolimits} }

\newcommand{\app}{ {\mathop{\mathrm {app}}\nolimits} }
\newcommand{\ext}{ {\mathop{\mathrm {ext}}\nolimits} }
\newcommand{\cur}{ {\mathop{\mathrm {cur}}\nolimits} }

\newcommand{\rel}{ {\mathop{\mathrm {rel}}\nolimits} }
\newcommand{\fnc}{ {\mathop{\mathrm {fnc}}\nolimits} }

\newcommand{\inter}{ {\mathop{\mathrm {int}}\nolimits} }
\newcommand{\cl}{ {\mathop{\mathrm {cl}}\nolimits} }
\newcommand{\lcl}{ {\mathop{\mathrm {lcl\,}}\nolimits} }
\newcommand{\rcl}{ {\mathop{\mathrm {rcl\,}}\nolimits} }
\newcommand{\cof}{ {\mathop{\mathrm {cof\,}}\nolimits} }
\newcommand{\add}{ {\mathop{\mathrm {add\,}}\nolimits} }
\newcommand{\sat}{ {\mathop{\mathrm {sat\,}}\nolimits} }
\newcommand{\tc}{ {\mathop{\mathrm {tc\,}}\nolimits} }
\newcommand{\unif}{ {\mathop{\mathrm {unif\,}}\nolimits} }

\newcommand{\olim}{\mathop{\mathrm{olim}}\nolimits}

\newcommand{\uhr}{\!\upharpoonright\!}
\newcommand{\lra}{ {\leftrightarrow} }
\newcommand{\ot}{ {\mathop{\mathrm {ot\,}}\nolimits} }
\newcommand{\ol}{\overline}
\newcommand{\cnc}{ {^\frown} }
\newcommand{\image}{\/``\,}

\newcommand{\scc}{\bm\upbeta}

\newcommand{\wh}{\widehat}
\newcommand{\wt}{\widetilde}
\newcommand{\inn}{\mathrm{in\,}}
\newcommand{\id}{\mathrm{id}}

\newcommand{\Lo}{\mathrm{Lo}}

\newcommand{\RC}{\mathrm{RC}}
\newcommand{\RCO}{\mathrm{RCO}}
\newcommand{\RCl}{\mathrm{RCl}}
\newcommand{\RO}{\mathrm{RO}}
\newcommand{\RClop}{\mathrm{RClop}}
\newcommand{\RegCl}{\mathrm{RegCl}}
\newcommand{\RegO}{\mathrm{RegO}}
\newcommand{\Cl}{\mathrm{Cl}}
\newcommand{\Clop}{\mathrm{Clop}}

\newcommand{\RK}{\mathrm{RK}}
\newcommand{\RF}{\mathrm{RF}}
\newcommand{\RB}{\mathrm{RB}}
\newcommand{\Comf}{\mathrm{C}}
\newcommand{\K}{\mathrm{K}}
\newcommand{\T}{\mathrm{T}}

\newcommand{\NCF}{\mathrm{NCF}}

\newcommand{\AC}{\mathrm{AC}}
\newcommand{\CH}{\mathrm{CH}}
\newcommand{\GCH}{\mathrm{GCH}}
\newcommand{\SCH}{\mathrm{SCH}}
\newcommand{\MA}{\mathrm{MA}}

\newcommand{\type}{\bm\uptau}

\author{Nikolai L.~Poliakov}%Nikolai L.~Poliakov}
%\author{}%Denis I.~Saveliev}

\title{On skew ultralimits and their applications in ultrafilter theory}
\address{HSE University}%HSE University}
\email{niknikols0@gmail.com}
%\address{%(corresponding author)
%Higher School of Modern Mathematics MIPT
% 1 Klimentovskiy per., Moscow, Russia}
%\address{%(corresponding author)
%}
%\email{d.i.saveliev@gmail.com}  
%\date{09.05.2026}
\maketitle
Ultrapowers and unions of chains are classic tools of model theory (see \cite{KC}). Combining these two constructions leads to new concepts, the most well-known of which is the operation of the \textit{ultralimit} (or \textit{limiting ultrapower}) of ordinal rank $\alpha$. In the sense of \cite{Kochen}, an ultralimit is the union of an increasing chain of successive ultrapowers and their unions with respect to a system of \textit{natural} (diagonal, at non-limit steps) embeddings. In this paper, we modify this notion by using a different (no less natural) system of embeddings. The resulting construction, the \textit{skewed ultralimit}, has useful applications in the theory of ultrafilters. Using it, \cite{PS2025Astana} obtained a generalization of Blass's theorem  on the model-theoretic characterization of the Rudin–Keisler order (see \cite{Blass}). In the present paper, we use skew ultralimits to give a precise description of the Rudin–Keisler (pre)order on the equivalence class of a Ramsey ultrafilter with respect to the Comfort preorder.
\par 
\section{Skew ultralimits and their basic properties}
For any $e:A\to B$ and ultrafilter $\mathfrak a\in \scc X$ define
$
e^\mathfrak a:
\prod_\mathfrak a A\to 
\prod_\mathfrak a B
$ 
by letting 
$$
e^\mathfrak a(g_\mathfrak a):=
(e\circ g)_\mathfrak a.$$
It is easy to see that the operation $e\mapsto e^\mathfrak a$ preserves such properties of $e$ as being the identity mapping, being an injection, or being a surjection. Moreover,
for all $e_0:A\to B$ and $e_1: B\to C$ we have 
$
e_1^\mathfrak a\circ e_0^\mathfrak a=(e_1\circ e_0)^\mathfrak a
$. Furthermore, for all models $\mathfrak A$ and $\mathfrak B$ if a mapping $e:\mathfrak A\to \mathfrak B$ is an embedding, an elementary embedding, an homomorphism, or an isomorphism, then so is the mapping  $e^\mathfrak a:
\prod_\mathfrak a\mathfrak A
\preceq 
\prod_\mathfrak a\mathfrak B$.
\begin{rem} For the diagonal embedding $d:\mathfrak M\to\prod_\mathfrak a\mathfrak M$, the embedding $d^\mathfrak a:\prod_a\mathfrak M\to \prod_\mathfrak a\prod_\mathfrak a\mathfrak M$ is not a diagonal embedding. 
\end{rem}
\begin{df}\label{df: unlim eng}
	For each model $\mathfrak M$, ultrafilter $\mathfrak a$, and ordinal $\alpha$, we define a direct system $(\mathcal M_\alpha, \mathcal F_\alpha)$ of elementary embeddings over the directed set $(\alpha+1, \leq)$ such that the direct systems $(\mathcal M_\beta, \mathcal F_\beta)$, $\beta< \alpha$, form an increasing chain:
	\begin{enumerate}
		\item[i.] $\mathcal M_0=\{\mathfrak M_0\}$ and $\mathcal F_0=\{e_{00}\}$, where $\mathfrak M_0=\mathfrak M$ and $e_{00}$ is the identity map.
		\item[ii.] $\mathcal M_1=\{\mathfrak M_0, \mathfrak M_1\}$, $\mathcal F_1=\{e_{00}, e_{01}, e_{11}\}$, where $\mathfrak M_0=\mathfrak M$, $\mathfrak M_1=\prod_\mathfrak a\mathfrak M$, $e_{00}$ and $e_{11}$ are the identity maps, and $e_{01}=d$ (the diagonal embedding),
		\item[iii.] if $\alpha$ is a limit ordinal, then $(\mathcal M_\alpha, \mathcal F_\alpha)$ is the direct limit of the union of the increasing chain $\{(\mathcal M_\beta, \mathcal F_\beta)\}_{\beta<\alpha}$,
		\item[iv.] if $\alpha=\beta+1\geq 2$, then $\mathcal M_\alpha=\mathcal M_\beta\cup \{\mathfrak M_\alpha\}$ and $\mathcal F_\alpha=\mathcal F_\beta\cup \{f_{\gamma\alpha}\}_{\gamma\leq\alpha}$, where
		\begin{enumerate}
			\item $\mathfrak M_\alpha=\prod_\mathfrak a\mathfrak M_\beta$,
			\item if $\gamma=\delta+1$, then $e_{\gamma\alpha}=e_{\delta\beta}^\mathfrak a$,
			\item if $\gamma=0$, then $e_{\gamma\alpha}= e_{1\alpha}\cdot e_{01}$,
			\item if $\gamma>0$ is a limit ordinal, then for any
			$
			\bm g\in M_{\mathfrak a,\gamma}
			$
			$$
			e_{\gamma\alpha}(\bm g)=
			e_{\delta\beta}^\mathfrak a(h_\mathfrak a)
			$$
			for some $\delta<\gamma$ and
			$
			h_\mathfrak a\in \bm g\cap M_{\mathfrak a,\delta+1},
			$
			\item $e_{\alpha\alpha}$ is the identity map.
		\end{enumerate}
	\end{enumerate}
	For any ordinal $\alpha$, the model $\mathfrak M_{\alpha}\in \mathcal M_\alpha$ will be called a \textit{skew ultralimit} (of $\mathfrak M$ of rank~$\alpha$ with respect to $\mathfrak a$), and denoted by $\mathfrak M_{\mathfrak a, \alpha}$.
\end{df}
\begin{thm}
	For any ordinal $\alpha$, ultrafilter $\mathfrak a$, and model $\mathfrak M$, the model $\mathfrak M_{\mathfrak a, \alpha}$ is well defined.
\end{thm}
\begin{rem}
	If in the definition of the skew ultralimit we replace all functions $e_{\alpha(\alpha+1)}$ with diagonal embeddings, we obtain the usual definition of the ultralimit. These two constructions have similar properties, but are not identical in general. A comparison of them can be found in \cite{Poliakov 2025}, where a general approach in category theory is also presented.
\end{rem}
\par The following theorem opens up the possibility of simple applications of the skew ultralimit.
\begin{thm}
	Let $\alpha$ be an ordinal, $\mathfrak a$ an ultrafilter, and $\{\mathfrak M_\beta\}_{\beta<\alpha}$ a sequence of models of the same signature $\sigma$ such that
	\begin{enumerate}
		\item $\mathfrak M_\beta\subseteq \mathfrak M_{\beta+1}$ for any ordinal $\beta<\alpha$,
		\item $\mathfrak N_\beta=\bigcup_{\gamma<\beta}\mathfrak M_\beta$  for any limit ordinal $\beta<\alpha$,
		\item there is a family $\{\iota_\beta\}_{\beta<\alpha}$ of isomorphisms $\iota_\beta:\prod_{\mathfrak a}\mathfrak M_\alpha\to \mathfrak M_{\alpha+1}$ such that all the diagrams
		$$
				\xymatrix{
			\prod_\mathfrak a\mathfrak M_0\ar[rr]^{\iota_0}&&\mathfrak M_1\\
			&&\\
		 \mathfrak M_0\ar[uu]_{\mathrm d}\ar[rruu]^{\mathrm id}&&
		}
		    	\xymatrix{
			&\\
			&\text{and}&\\
		    	 &&
	    	 }
	        	\xymatrix{
			\prod_\mathfrak a\mathfrak M_{\beta+1}\ar[rr]^{\iota_{\beta+1}}&&\mathfrak M_{\beta+2}\\
			&&\\
			\prod_\mathfrak a\mathfrak M_\beta \ar[uu]_{\mathrm id}\ar[rr]^{\iota_\beta} &&\mathfrak M_{\beta+1}\ar[uu]^{\mathrm id	}
		}
		$$
		are commutative, where $\mathrm d$ is the diagonal embedding, and $\mathrm id$ identity mapping.
	\end{enumerate}
Then $\mathfrak M_\beta\cong(\mathfrak M_0)_{\mathfrak a,\beta}$ for all $\beta<\alpha$.
\end{thm}
\section{Applications to the theory of ultrafilters}
Recall that the set of ultrafilters $\scc X$ on $X$ endowed with a natural topology with an clopen base of the sets $\{\mathfrak u\in \scc X: A\in \mathfrak u\}$, $A\subseteq X$, is the \textit{\v{C}ech--Stone compactification} of the discrete space $X$. This means that $X$~is dense in~$\scc X$, and every (trivially continuous) map~$h$ of~$X$ into any compact Hausdorff space~$Y$ uniquely extends 
to a~continuous map~$\wt h$ of~$\scc X$ into~$Y$:
	$$
\!\!\!\!\!\!\!\!\!\!\!\!\!\!\!\!\!\!\!\!
\xymatrix{
	&\scc X\,
	\ar@{-->}^{\wt{h}\quad}[drr]&&
	\\
	&X\,
	\ar[u]
	\ar[rr]^{h}
	&&\,Y&
}
$$ 
\par Classical objects of theory of ultrafiltars are various (pre)orders on $\scc X$.
\par For all $\mathfrak u,\mathfrak v\in \scc X$, $\mathfrak u$ is \textit{Rudin--Keisler less} 
than~$\mathfrak v$, 
denoted by $\mathfrak u\le_\RK\mathfrak v$,
iff there exists 
$f:\omega\to\omega$ such that 
$\widetilde f(\mathfrak v)=\mathfrak u$. Equivalence classes $\tau(\mathfrak u)$, $\mathfrak u\in \scc X$, with respect to the relation ${\approx}_\RK:={\le}_\RK\cap {\le}^{-1}_\RK$ are called \textit{ultrafilter types}. Therefore, $\le_\RK$ generates 
the \textit{Rudin--Keisler order}
on the set of types of ultrafilters. A non-principal ultrafilter $\mathfrak u\in \scc X$ is called ($\RK$-)\textit{minimal} if $\mathfrak v\le_\RK \mathfrak u$ implies $\mathfrak v\approx_\RK \mathfrak u$ for any non-principal ultrafilter $\mathfrak v\in \scc X$. Minimal ultrafilters have numerous alternative characterizations, see \cite{CN} and \cite{Poliakov Rams uf}. In particular, they are known as \textit{Ramsey ultrafilters}.
\par The Comfort pre-order is defined as follows (for simplicity, we restrict ourselves to the case of ultrafilters on $\omega$). For $\mathfrak u,\mathfrak v\in\scc \omega$ 
\begin{enumerate}
	\item[(i)] 
	a~space~$Y$ is \textit{$\mathfrak u$-compact} iff 
	$\widetilde f(\mathfrak u)\in Y$ for any $f:\omega\to Y$;
	\item[(ii)] 
	$\mathfrak u$ is \textit{Comfort less} 
	than $\mathfrak v$, 
	denoted by $\mathfrak u\le_\mathrm C\mathfrak v$, iff any $\mathfrak v$-compact space 
	is $\mathfrak u$-compact.
\end{enumerate} 
   Equivalence classe of ultrafilter $\mathfrak u\in \scc\omega$ with respect to the relation ${\approx}_{\mathrm C}:={\le}_{\mathrm C}\cap {\le}^{-1}_{\mathrm C}$ is denoted by $T_\mathrm C(\mathfrak u)$. We also use the notation $T_\mathrm C^\tau(\mathfrak u):=\{\tau(\mathfrak u): \mathfrak u\in T_\mathrm C(\mathfrak u)\}$.
   \par Numerous papers are devoted to the study of these (and other) preorders on sets $\scc X$ (and in particular, on the set $\scc\omega$). A survey and bibliography (far from complete) can be found in~\cite{PS2026_1} and~\cite{PS2026_2}. The connection between the Comfort and Rudin-Keisler preorders is deeply studied in~\cite{Ferreira,Ferreira1} (see also \cite{PS2025Astana}). We also note one result of the paper \cite{Booth}, which is directly related to our work: $(\scc\omega, \le_\RK)$ contains a fragment isomorphic to the ultrapower of $(\omega, \leq)$. Note that most studies establish the embeddability of certain orders in $(\scc, \le_\RK)$ or a set $\scc X$ preordered by some other relation. It is very rare to be able to accurately describe the Rudin--Keisler preorder on some meaningful set of ultrafilters. Our study provides just such an example.
   \par Let $\mathrm{DF}$ be a~set of all injective functions 
   $f:\omega\to \scc\omega$ with a~discrete range.
   \begin{df}\label{df: W}
   	For any ultrafilter $\mathfrak u\in \scc \omega\setminus\omega$ and ordinal $\alpha$ define the sets $W_{<\alpha}(\mathfrak u)$ and $W_\alpha(\mathfrak u)$ by recursion on $\alpha$:
   	\begin{enumerate}
   		\item $W_0(\mathfrak u)=\omega$, $W_{<0}(\mathfrak u)=\emptyset$,
   		\item if $\alpha>0$, $W_{<\alpha}(\mathfrak u)=\bigcup_{\beta<\alpha}W_\beta(\mathfrak{u})$ and $W_\alpha(\mathfrak{u})$ is the set of all ultrafiters $\widetilde f(\mathfrak u)$ where $f\in DF$ and for all $\beta<\alpha$ $$\{i<\omega: f(i)\in W_{<\alpha}(\mathfrak{u})\setminus W_{<\beta}(\mathfrak{u})\}\in \mathfrak u.$$
   	\end{enumerate}
   	Denote $W(\mathfrak{u}):=\bigcup\limits_{\alpha\in \mathrm{Ord}}W(\mathfrak u)$, $W^\tau_{<\alpha}(\mathfrak{u}):=\{\tau(\mathfrak u): \mathfrak u\in W_{<\alpha}(\mathfrak{u})\}$ and $W^\tau(\mathfrak{u}):=\{\tau(\mathfrak u): \mathfrak u\in W(\mathfrak{u})\}$.
   \end{df}
This construction is closely related to relations $R_\alpha$ introduced in \cite{PS2025Astana}, as well as to the concept of \textit{tower} of an ultrafilter $\mathfrak u$ from \cite{Ferreira}.
\begin{thm}
	For any non-principal ultrafilter $\mathfrak u\in \scc\omega$
	\begin{enumerate}
		\item $W(\mathfrak{u})=W_{<\omega_1}(\mathfrak{u})\subseteq T_\mathrm{C}(\mathfrak u)$,
		\item if $\mathfrak u$ is a Ramsey ultrafilter, then $W^\tau(\mathfrak{u})= T^\tau_\mathrm{C}(\mathfrak u)$.
	\end{enumerate}
\end{thm}
\begin{thm}\label{th: W and M}
	For any non-principal ultrafilter $\mathfrak u\in \scc\omega$ and ordinal $\alpha$, the poset $(W^\tau_{<\omega+\alpha}(\mathfrak{u}), \le_\RK)$ is isomorphic to the skew ultralimit $\mathfrak M_{\mathfrak u, \alpha}$ of the poset $\mathfrak M=(\omega, \leq)$ with the standard order.
\end{thm}
It is easy to see that for every non-principal ultrafilter $\mathfrak u\in \scc\omega$ and model of countable signature~$\sigma$, the model $\mathfrak M_{\mathfrak u, \omega_1}$ is $\omega_1$-saturated. Therefore, the following fact holds.
\begin{cor}\label{th: W and M omega 1}
	For any non-principal ultrafilter $\mathfrak u\in \scc\omega$
	$$
	(W^\tau(\mathfrak{u}), \le_\RK)\cong \prod_\mathfrak u(\omega, \leq).
	$$
	If $\mathfrak u$ is a Ramsey ultrafilter, then $$(T^\tau_\mathrm{C}(\mathfrak{u}), \le_\RK)\cong \prod_\mathfrak u(\omega, \leq).$$
\end{cor}
\begin{rem}
	Previously, in \cite{P_Erl_2024}, the Rudin-Keisler preorder on the $\approx_\mathrm C$-equivalence class of the Ramsey ultrafilter $\mathfrak u\in \scc\omega$ was characterized using the concept of the $\mathrm{o}$-limit.
\end{rem}

\end{document}